\documentclass[10pt, a4paper]{article}

\usepackage[utf8]{inputenc}
\usepackage{amsmath} 
\usepackage{amssymb} 
\usepackage{amsthm}
\usepackage[T1]{fontenc}
\usepackage[french]{babel}

\usepackage[backend=bibtex,firstinits=true, isbn=false, doi=false]{biblatex}
\renewbibmacro{in:}{}
\newtheorem{nota}{Notation} 

\newtheorem*{aff}{Affirmation}
\newtheorem{proposition}{Proposition}[section]
\newtheorem{Remarque}[proposition]{Remarque}

\newtheorem{fait}[proposition]{Fait}
\newtheorem{definition}[proposition]{Définition}
\newtheorem{theoreme}[proposition]{Théorème}
\newtheorem{Lemme}[proposition]{Lemme}
\newtheorem{corollaire}[proposition]{Corollaire}

\newcommand{\ad}{\operatorname{ad}}
\newcommand{\Int}{\operatorname{Int}}
\newcommand{\engendre}[1]{\left\langle#1\right\rangle}

\title{Quelques résultats sur les anneaux de Lie qui n'ont pas de chaîne infinie de centralisateurs\footnote{ mots-clé : algèbre de Lie, radical, condition de chaîne sur les centralisateurs
\\
MSC 2020 : 20N05, 20F11, 03C60}}

\author{Samuel Zamour}
\date{18 mars 2025}

\begin{document}
\maketitle
\begin{abstract}
   Dans la lignée des résultats existants pour les groupes, on montre que le radical nilpotent existe dans les anneaux de Lie qui n'ont pas de chaîne infinie de centralisateurs. On établit également un analogue du théorème d'Engel si la caractéristique est nulle.
\end{abstract}

\section{Introduction}
Les algèbres de Lie sont l'objet d'un regain d'attention de la part des théoriciens des modèles, comme en témoignent les deux articles récents \cite{DT1} et \cite{DT2} portant sur les algèbres de Lie rangées.  Il s'agit de prolonger cette étude et de généraliser au contexte des anneaux de Lie certains résultats connus pour les groupes satisfaisant de bonnes propriétés du point de vue de la théorie des modèles. On s'intéressera aux anneaux de Lie qui n'ont pas de chaîne infinie de centralisateurs. 
\begin{definition}
\begin{itemize}
\item Un anneau de Lie $(L,+, [ , ])$ est un groupe abélien $(L,+)$ muni d'une opération bi-additive $[ , ]$ qui vérifie les propriétés suivantes :
\begin{itemize}
    \item Pour tout $x\in L$, $[x,x]=0$.
    \item (identité de Jacobi) Pour tous $x,y,z\in L$, \[[x,[y,z]]+[z,[x,y]]+[y,[z,x]]=0.\]
\end{itemize}
\item Si on ajoute une structure d'espace vectoriel sur $(L,+)$ de sorte que $[,]$ soit bilinéaire, on parle d'algèbre de Lie. La caractéristique de $L$ correspond à celle du corps sous-jacent.
\item Soit $X\subset L$, on définit le centralisateur de $X$ par $C_L(X)=\{y\in L : \forall x\in X, [y,x]=0\}$. On dit que $L$ satisfait la condition de chaîne minimale sur les centralisateurs s'il n'existe pas de chaîne infinie descendante ou ascendante de centralisateurs. On parlera d'anneau de Lie $\mathfrak{M}_c$ .
\end{itemize}
\end{definition}

Les groupes $\mathfrak{M}_c$ ont été étudiés par Bryant dans \cite{Bry}; l'investigation de leurs propriétés a été poursuivie par Wagner dans \cite{Wag1}, \cite{Wag2} et \cite{DW}. Les algèbres de Lie vérifiant certaines conditions de chaîne ont déjà été considérées (voir \cite{St}); mais à notre connaissance, la condition de chaîne sur les centralisateurs n'a jamais été étudiée dans ce contexte.

Dans une première section, on généralise au contexte des anneaux de Lie certains résultats sur les centralisateurs itérés. Dans une deuxième section, on montre que le radical nilpotent existe dans les anneaux de Lie $\mathfrak{M}_c$. Dans la troisième section, on caractérise les éléments ad-nilpotents dans une algèbre de Lie $\mathfrak{M}_c$ de caractéristique nulle et on établit dans ce contexte un analogue du théorème d'Engel. Enfin, dans une quatrième et dernière section, on donne quelques résultats relativement à la définissabilité des radicaux dans le cas stable.
\section{Préliminaires sur les anneaux et les algèbres de Lie}
Pour les notions de base sur les anneaux de Lie et les algèbres de Lie, on pourra se reporter aux premiers chapitres de \cite{Hum} ou bien à l'introduction de \cite{DT1}. On adopte le language des anneaux de Lie, mais tous les résultats de cette section se transfèrent de façon immédiate aux algèbres de Lie.
Commençons par quelques définitions et propriétés élémentaires pour fixer notre cadre.
\begin{definition} Soit $L$ un anneau de Lie.
    \begin{enumerate}
        \item Un sous-groupe $H$ est un sous-anneau si pour tous $x,y\in H$, $[x,y]\in H$. C'est un idéal si de plus, $[x,H]\subseteq H$ pour tout $x\in L$.
        \item Pour $X\subset L$, on définit $[X,X]=X'$ comme le sous-groupe additif engendré par les commutateurs $[x,y]$. On définit par récurrence la suite centrale descendante et la suite dérivée :
        $L^1=L'$, $L^{n+1}=[L,L^n]$ et $L^{(1)}=L'$, $L^{(n+1)}=[L^{(n)},L^{(n)}]$.
        \item On définit par récurrence les centralisateurs itérés de $A/H$, où $A$ est une sous-anneau et $H$ un idéal de $A$ : 
        \[ C_L^0(A/H)=H \]
        \[C_L^{n+1}(A/H)=\{x\in \bigcap_{1\leq i\leq n}N_L(C_L^i(A/H) : [x,A]\leq C_L^n(A/H)\}.\]
        Dans le cas où $A=L$ et $H=\{0\}$, on parlera de centres itérés, notés $Z_i(L)$.
    \end{enumerate}
\end{definition}
\begin{nota}
    On utile le symbole $\subseteq$ pour l'inclusion ensembliste, $\leq_{gr}$ pour les sous-espaces ou les sous-groupes, $\leq $ pour les sous-anneaux ou les sous-algèbres et $\lhd$ pour les idéaux. Par ailleurs, étant donné un sous-ensemble $X$, on définit $\engendre{X}$ comme le sous-anneau engendré par $X$.
\end{nota}
La notion d'idéal constitue un analogue de celle de sous-groupe normal, mais les homomorphimes de groupes $\ad_x$, définis par $\ad_x(y)=[x,y]$, ne sont pas des automorphismes. En revanche, en raison de l'identité de Jacobi, ce sont des dérivations; on rappelle qu'un homomorphisme additif $D$ est une dérivation d'un anneau de Lie si $D([x,y])=[D(x), y ]+[x,D(y)]$. 
\begin{definition}
    Soit $L$ un anneau de Lie. On dit qu'un élément $x$ est $n$-ad-nilpotent si $\ad_x^n(L)=0$ pour un certain entier $n$. L'anneau de Lie $L$ est ad-nilpotent si tous ses éléments sont ad-nilpotents.
\end{definition}

Il existe également une notion d'idéal caractéristique :
\begin{definition}
    Soit $I\lhd L$ un idéal. On dit qu'il est caractéristique s'il est stable sous toute dérivation de $L$.
\end{definition}
    \begin{fait}
        Les sous-groupes $L^{(n)}, L^n, Z_n(L)$ sont des idéaux caractéristiques.
    \end{fait}
Nous allons démontrer certains résultats sur les centralisateurs itérés. Il s'agit de généralisations immédiates de ceux connus pour les groupes; les démonstrations sont parfois un peu plus simples en raison de l'identité de Jacobi.
\begin{definition}
Soit $L$ un anneau de Lie. Pour un sous-groupe $A$, $[A,_n x]$ désigne le sous-groupe engendré par les éléments de la forme $\ad_{a_1}\circ\dots\circ\ad_{a_n}(x) $ pour $a_1,\dots,a_n\in A$.
\end{definition}
\begin{proposition}\label{centralisateur itéré}
Les centralisateurs itérés sont des sous-anneaux. De plus, les centralisateurs itérés de la forme $C_L^n(I)$, où $I$ est un idéal, sont des idéaux et $C_L^n(I)=\{ x\in L : [I,_n x]=0\}$.
    
\end{proposition}
\begin{proof}
    On procède par récurrence.
    Pour $n=1$, $C_L(A/H)=\{x\in N_L(H) : [x,A]\leq_{gr} H\}$, c'est un sous-anneau car pour tous $x,y\in C_L(A/H)$ et $z\in A$, on a : $[x+y,z]=[x,z]+[y,z]\in H$ et $[[x,y],z]=([[x,z],y]+[x,[y,z]])\in H$, car $H$ est un sous-groupe et $y,x\in N_L(H)$.
    Supposons la propriété vérifiée pour $n$. Comme les normalisateurs sont des sous-anneaux, par hypothèse de récurrence, les mêmes calculs nous donnent que $C_L^{n+1}(A/H)$ est un sous-anneau. On procède de même pour montrer que les centralisateurs de la forme $C_L^n(I)$, où $I$ est un idéal, sont des idéaux

    On montre par récurrence que  $C_L^n(I)=\{ x\in L : [I,_n x]=0\}$. C'est évident pour $n=1$. Supposons la propriété démontrée pour $n$. Comme les centralisateurs itérés sont des idéaux, $x\in C_L^{n+1}(I)$ ssi $[x,I]\leq_{gr} C_L^n(I)$. Si $[I,_{n+1} x]=0$, alors pour tous $y_1,\dots,y_{n},y_{n+1}\in I$, $(\ad_{y_1}\circ...\circ \ad_{y_n})(\ad_{y_{n+1}}(x))=0$ et donc par hypothèse de récurrence $\ad_{y_{n+1}}(x)\in C_L^n(I)$. Réciproquement, si $[I,x]\leq_{gr} C_L^n(I)$, alors par hypothèse de récurrence, $[I,_{n+1} x]=0$.
\end{proof}
\begin{Remarque}
    Pour un sous-anneau $A$, on a que $C_L^n(A)\subseteq \{x\in L : [A,_n x]=0\}$.
\end{Remarque}
\begin{Lemme}
Soit $L$ un anneau de Lie, et soit $A$ un sous-anneau. Alors $N_L(A)\leq N_L(C_L^n(A))$, pour tout entier $n$.
\end{Lemme}
\begin{proof}
On procède par récurrence. Pour $n=1$, soient $y\in N_L(A)$, $x\in C_L(A)$ et $a\in A$; alors $[[y,x],a]=[[y,a],x]+[y,[x,a]]=0$ car $[y,a]\in A$ et $[x,a]=0$. 

Supposons la propriété vérifiée pour $n$. Soient $y\in N_L(A)$ et $x\in C_L^{n+1}(A)$; alors $x\in \bigcap_{i\leq n} N_L(C^i(A))$, mais par hypothèse de récurrence, $y$ est aussi contenu dans $\bigcap_{i\leq n} N_L(C^i(A))$, et ainsi $[y,x]\in \bigcap_{i\leq n} N_L(C^i(A))$. De plus, pour $a\in A$, le même calcul qu'au paragraphe précédent et l'hypothèse de récurrence nous donne que $[[y,x],a]\in C_L^n(A)$.
\end{proof}
En raison de l'identité de Jacobi, l'équivalent du Lemme des Trois Groupes est immédiat.
\begin{fait}
Si $X,Y$ et $Z$ sont des sous-ensembles de $L$  tels que pour un sous-groupe $H$, $[[X,Y],Z]\leq_{gr} H$ et $[[Y,Z],X]\leq_{gr} H$, alors $[[Z,X],Y]\leq_{gr} H$.
\end{fait}
Pour les deux Lemmes suivants, nous reprenons presque littéralement les démonstrations données par Wagner dans \cite{Wag1} dans le cas des groupes.
\begin{Lemme}[{\cite[Lemma 2.4]{Bry}, \cite[Lemma 0.1.11]{Wag1}}]\label{cent it 1}  
   Soit $H$ un sous-anneau de $L$. Alors pour tout $i<j$ :
   \[ [H^i, C_L^j(H)]\leq_{gr} C_L^{j-i-1}(H).\]
   En particulier, $H^{i-1}$ et $C_L^i(H)$ commutent.
\end{Lemme}
\begin{proof}
On procède par récurrence sur $i+j$. Pour $i+j=1$, on a $i=0$ et $j=1$, le lemme est vrai. Supposons que la propriété est vérifiée pour $i+j=n$, et considérons $i,j$ tels que $i+j=n+1$. Alors, par récurrence :
\[ [[H^{i-1},C_L^{j}(H)], H]\leq_{gr} [C_L^{j-i},H]\leq_{gr} C_L^{j-i-1}(H)\]
 et également,
 \[ [[H, C_L^{j}(H)], H^{i-1}]\leq_{gr} [C_L^{j-1}(H), H^{i-1}]\leq_{gr} C_L^{j-i-1}(H)\]
 Par conséquent, en vertu de l'identité de Jacobi, on a :
 \[[H^i, C_L^j(H)]=[[H^{i-1},H],C_L^{j}(L)]\leq_{gr} C_L^{j-i-1}(H).\qedhere\]
\end{proof}
\begin{Lemme}[{\cite[Lemma 2.5]{Bry}, \cite[Lemma 0.1.12]{Wag1}}]\label{cent it 2} 
  Supposons que $K\leq H \leq L$, et $C_L(K^i)=C_L(H^i)$ pour tout $i<j$. Alors $C_L^j(K)=C_L^j(H)$.  
\end{Lemme}
\begin{proof}
Par récurrence sur $j$. C'est évident pour $j=0,1$. Supposons que la propriété est vraie pour $j-1$.
\begin{aff}
Pour tout $i<j$, on $[H^{j-i-1},C_L^j(K)]\leq C_L^{i}(H)$.  
\end{aff} 
\begin{proof}
Par récurrence sur $i$. Pour $i=0$, on doit montrer que $[H^{j-1},C_L^j(K)]=0$, ce qui est vrai par le Lemme \ref{cent it 1}. Supposons que l'affirmation est vraie pour $i-1$. On a 
\[[[H^{j-i-1},K],C_L^j(K)]\leq_{gr} [H^{j-i},C_L^j(K)]\leq_{gr} C_L^{i-1}(H).\]
De plus, par le Lemme \ref{cent it 1} et par récurrence, on obtient : 
\[[[K, C_L^j(K)], H^{j-i-1}]\leq_{gr} [C_L^{j-1}(K), H^{j-i-1}]=[C_L^{j-1}(H),H^{j-i-1}]\leq_{gr} C_L^{i-1}(H).\] 
D'après l'identité de Jacobi,
\[[[H^{j-i-1}, C_L^j(K)],K]\leq_{gr} C_L^{i-1}(K).\]
Mais $H^{j-i-1}$ et $C_L^j(K)$ normalisent $C_L^k(K)=C_L^k(H)$ pour tout $k<i$, par conséquent 
\[[H^{j-i-1},C_L^j(K)]\leq_{gr} C_L^i(K)=C_L^i(H).\]
\end{proof}
En particulier, $[H, C_L^j(K)]\leq_{gr} C_L^{j-1}(H)$, et $C_L^j(K)$ normalise $C_L^i(H)=C_L^i(K)$ pour tout $i<j$. Ainsi, $C_L^j(K)\leq C_L^j(H)$. D'autre part, $[K,C_L^j(H)]\leq_{gr} C_L^{j-1}(H)=C_L^{j-1}(K)$, et $C_L^j(H)$ normalise $C_L^i(K)=C_L^i(H)$ pour tout $i<j$. Par conséquent, $C_L^j(H)\leq C_L^j(K)$.
\end{proof}

\section{Radical Nilpotent}
\begin{definition}
    On dit que $L$ est $\mathfrak{Z}_f$ si pour chaque centre itéré $Z_i(L)$, le centre de $L/Z_i(L)$ est le centralisateur d'un nombre fini d'éléments.
\end{definition}
On reprend tels quels les arguments de Wagner pour établir :
\begin{theoreme} [{\cite[Theorem 1.2.1]{Wag1}}]
    Un anneau de Lie $\mathfrak{M}_c$ est $\mathfrak{Z}_f$.
\end{theoreme}
\begin{proof}
    Pour chaque $n$, la condition de chaine sur les centralisateurs nous donne l'existence d'un ensemble fini $\overline{x}_n\in L^n$ tel que $C_L(L^n)=C_L(\overline{x}_n)$. De plus, il existe un sous-ensemble fini $\overline{k}_n\in L$ tel que $\overline{x}_n\in \engendre{\overline{k}_n}^n$. On considère le sous-anneau $K=\engendre{\overline{k}_0,...,\overline{k}_n}$; alors $K\leq L$, et pour $i\leq n$, on a $\overline{x}_i\in K^i$, i.e., $C_L(K^i)=C_L(L^i)$. Donc, d'après le Lemme \ref{cent it 2}, $C_L^n(K)=C_L^n(L)=Z_n(L)$. Pour $\overline{L}=L/Z_{n-1}(L)$, on a $Z(\overline{L})=C_{\overline{L}}(K)=C_{\overline{L}}(\overline{k})$, où $\overline{k}$ est l'image de l'ensemble engendrant $K$. 
\end{proof}
On énonce une caractérisation des anneaux de Lie $\mathfrak{M}_c$ localement nilpotents parfaitement analogue à celle valable pour les groupes.
\begin{theoreme}[{\cite[Corollary 2.2]{Bry}}]\label{localement nilpotent} 
Soit $L$ un anneau de Lie $\mathfrak{M}_c$ localement nilpotent. Alors $L$ est résoluble.
\end{theoreme}
\begin{proof}
On peut supposer que tout centralisateur de la forme $C_L(x)$ qui est un sous-anneau propre est résoluble. On peut également supposer que $L$ n'est pas abélien.

Il existe $x\in Z_2(L)-Z(L)$ : en effet, $Z_2(L)/Z(L)=Z(L/Z(L))=Z(A)=C_A(a_1,...,a_n)$ car $L$ est $\mathfrak{Z}_f$, mais $A$ est localement nilpotent, et donc $Z(\engendre{a_1,...,a_n})\neq \{0\}$. L'homomorphisme $\ad_x$ définit un homomorphisme d'anneau de $L$ vers $Z(L)$ : $[x,[y,z]]=[[x,y],z]+[y,[x,z]]=0=[[x,y],[x,z]]$, car $[x,y],[x,z]\in Z(L)$. Par conséquent, $C_L(x)$ est un idéal propre de $L$ et donc par induction, c'est un idéal résoluble. On a $L/C_L(x)\simeq Z(L)$, et on obtient finalement que $L$ est résoluble.
\end{proof}
\begin{definition}
    On définit l'idéal de Fitting par :
    \[F(L)=\sum \{I : \text{I est un idéal nilpotent}\}.\]
\end{definition}
\begin{Lemme}
    $F(L)$ est localement nilpotent.
\end{Lemme}
\begin{proof}
    Soient $x_1, \dots , x_n\in F(L)$, le sous-anneau $\engendre{x_1,\dots x_n}$ est contenue dans $I_1+...+I_n$, où $I_1,\dots,I_n$ sont des idéaux nilpotents; mais $I_1+\dots+I_n$ est un idéal nilpotent, d'après \cite[Lemma 0.2]{St}.
\end{proof}
C'est le point novateur de cette section; on adapte à notre contexte les arguments de \cite{DW} qui utilisent le formalisme des anneaux de groupes.
\begin{definition}
    \begin{itemize}
        \item On définit une représentation d'un anneau de Lie $L$ sur un groupe abélien $V$, comme la donnée d'un homomorphisme d'anneaux de Lie  $\phi : L \rightarrow End(V) $, où le crochet de $End(V)$ est donné par $[f,g]=f\circ g - g\circ f$. On dit que $V$ est un $L$-module.
        \item  Soit $V$ un $L$-module. On définit une structure d'anneau de Lie sur le produit $L\times V$, on parlera du produit semi-direct $L\rtimes V$, par l'opération : $[(g,v),(g',v')]=([g,g'], \phi(g)(v')-\phi(g')(v))$.
    \end{itemize}
\end{definition}
\begin{Lemme}\label{produit semi direct} [{à comparer avec \cite[Lemma 3]{DW}}]

    Soit $\phi : L\rightarrow End(V)$ la donnée d'un $L$-module, où $L$ est abélien. On note $A=L\rtimes V$ le produit semi-direct d'anneaux de Lie correspondant. On suppose que $C_V(L)=C_V(x_1,...,x_k)$. Soit $v\in V$ tel que pour tout $1\leq i \leq k$, il existe un entier $n_i$ avec $\ad^{n_i}_{x_i}(v)=0$. Alors $[L,_m v]=0$ pour $m=1+\sum_{i=1}^{k} (n_i-1)$.
\end{Lemme}
\begin{proof}
On prend $m=1+\sum_{i=1}^{k} (n_i-1)$. En tant que que module de Lie, la condition sur les commutateurs se traduit par $\phi(x_i)^{n_i}(v)=0$. Par conséquent, pour tout choix de $m$ indices $i_1,...,i_m$ parmi $\{1,\dots,k\}$, au moins un indice doit apparaitre $n_i$ fois, mais comme $L$ est abélien, les endomorphismes $\phi(x_i)$ commutent entre eux. Ainsi,
\[\phi(x_{i_1})\circ\dots\circ \phi(x_{i_m})(v)=0.\]
Par conséquent, $\phi(x_i)\circ\phi(x_{i_2})\circ \dots \circ \phi(x_{i_m})(v)=0$ pour tout $i\in \{1,...,k\}$. Donc pour tout $x'_1\in L$, on obtient : \[ \phi(x_1')\circ \phi(x_{i_2})\circ \dots \phi(x_{i_m})(v)=0.\]
En utilisant la commutativité du module de Lie et en répétant $m$  fois l'opération, on obtient pour tous éléments $x'_1,...,x'_m\in L$ :
\[\phi(x'_1)\circ\dots\circ \phi(x'_m)(v)=0.\]
Par conséquent, $[L,_m v]=0$.  
\end{proof}
\begin{proposition}[{\cite[Proposition 6]{DW}}]\label{critère nilpotence}
  Soit $L$ un anneau de Lie $\mathfrak{M}_c$ résoluble tel que pour tout $x\in L$, il existe un entier $n$ tel que $\ad_x^n(y)=0$ pour tout $y\in L'$. Alors $L$ est nilpotent.  
\end{proposition}
\begin{proof}
    On procède par récurrence sur la classe de résolubilité, le cas abélien étant trivial. On peut donc supposer que $L'$ est nilpotent.
    Il suffit de montrer par récurrence que tout centre itéré de $L'$ est contenu dans un centre itéré de $L$ : il existera alors un indice $i$ tel que $L'=Z_n(L')\leq Z_i(L)$ et donc $L/Z_i(L)$ est abélien et $L=Z_{i+1}(L)$.
    Evidemment, $\{0\}=Z_0(L')$ est contenu dans $Z_0(L)$. Supposons que $Z_i(L')\leq Z_k(L)$; les sous-anneaux $Z_{i+1}(L')$ et $Z_i(L')$ sont des idéaux de $L$. On considère $A=Z_{i+1}(L')/Z_i(L')$; c'est un anneau de Lie abélien sur lequel $L$ agit de façon adjointe, de sorte que $A$ est un $L$-module. L'anneau de Lie $A$ est centralisé par $L'$ et on obtient donc une structure de $L/L'$-module.
    
    Comme $L$ est $\mathfrak{Z}_f$ et $A\leq L/Z_i(L)$, il existe $x_1,\dots,x_n$ dans $L$ tels que $C_A(L)=C_A(x_1,..,x_n)$. Pour tout $x_i$, il existe un entier $n_i$ tel que $\ad_{x_i}^{n_i}(y)=0$ pour tout $y\in L'$. De plus, pour tout $y\in L'$ et $a\in A$, on a bien $\ad_{x_i+y}(a)=0$ ssi $\ad_{x_i}(a)=0$, et donc $C_A(L/L')=C_A(\overline{x_1},\dots,\overline{x_n})$.
    
    D'après le Lemme \ref{produit semi direct}, pour $m=1+\sum_{i=1}^{k} (n_i-1)$, on obtient pour tout $a\in A$, $[L/L',_ma]=[L,_m a]=0$; par conséquent, pour tout  $y\in Z_{i+1}(L')$, on a $[L,_m y]\leq Z_i(L')\leq Z_k(L)$. Par conséquent, d'après la Proposition \ref{centralisateur itéré}, $y\in Z_{k+m}(L)$. Finalement, $Z_{i+1}(L')\leq Z_{k+m}(L)$.
\end{proof}
\begin{corollaire}\label{ad-nilpotence résoluble}
     Soit $L$ un anneau de Lie $\mathfrak{M}_c$ résoluble et ad-nilpotent. Alors $L$ est nilpotent.
 \end{corollaire}
 On en déduit également la nilpotence de l'idéal de Fitting.
\begin{theoreme}[{\cite[Theorem 8]{DW}}]
    Soit $L$ un anneau de Lie $\mathfrak{M}_c$. Alors $F(L)$ est un idéal nilpotent.
\end{theoreme}
\begin{proof}
    L'idéal $F=F(L)$ est localement nilpotent et donc résoluble d'après la Proposition \ref{localement nilpotent}. Soit $x\in F$; cet élément appartient à un idéal nilpotent $I$ de $L$, de classe de nilpotence $r$ disons. Par conséquent, pour un élément $y$ de $F'$, $[x,y]\in I$ et donc $ad^r_{x}(y)=0$. On conclut d'après la Proposition \ref{critère nilpotence}.
\end{proof}
\section{Théorème de Engel}
Dans cette section, nous allons établir un analogue du théorème de Engel. Sauf mention explicite du contraire, $L$ est une algèbre de Lie sur un corps de caractéristique nulle.
D'après les travaux de Hartley, il est possible de définir un analogue du radical de Baer dans ce contexte. 
\begin{definition}
    Soit $L$ une algèbre de Lie. On dit que $H$ est un sous-idéal de $L$ s'il existe une suite finie de sous-algèbres telle que :
    \[H_0=H\lhd H_1 \lhd\dots\lhd H_n=L.\]
    L'indice d'un sous-idéal désigne la longueur minimale d'une suite de la forme précédente.
\end{definition}
\begin{fait} \cite{Hart} Le radical de Baer $B(L)$ est la sous-algèbre engendrée par les sous-idéaux nilpotents et de dimension finie. On a :
\begin{enumerate}
    \item\
    
    \begin{tabular}{rcl}
    $B(L)$ & $ =$ & la sous-algèbre engendrée par les sous-idéaux nilpotents\\
    & $=$ & la sous-algèbre engendrée par les sous-idéaux de dimension 1\\
    & $=$ & l'ensemble des x dans $L$ tels que $<x>$ est un sous-idéal.
    \end{tabular}
    \item $B(L)$ est un idéal localement nilpotent caractéristique.
    \item $F(L)$ est un idéal caractéristique.
\end{enumerate}
    \end{fait}
\begin{Lemme}[{\cite[Corollary 10]{DW}}]\label{sous-idéal nilpotent}
Soit $L$ une $\mathfrak{M}_c$ algèbre de Lie de caractéristique nulle. Alors toute famille d'idéaux nilpotents engendre un idéal nilpotent. Un sous-idéal nilpotent est contenu dans un idéal nilpotent.
\end{Lemme}
\begin{proof}
Toute famille d'idéaux nilpotents est contenue dans l'idéal de Fitting qui est nilpotent, donc engendre nécessairement un idéal nilpotent.

Pour la deuxième assertion, on procède par récurrence sur la classe du sous-idéal; soit $I=I_0\lhd \dots\lhd I_n=L$. Si $n=0$, $I$ est nilpotent. Supposons que l'affirmation est vraie pour $n-1$. Alors $I$ est un sous-idéal de $I_{n-1}$ et par récurrence il est contenu dans $F(I_{n-1})$ qui est un idéal caractéristique dans $I_{n-1}$; c'est donc un idéal nilpotent de $L$.
\end{proof}
\begin{corollaire}
    Soit $L$ une $\mathfrak{M}_c$ algèbre de Lie de caractéristique nulle. Alors $B(L)=F(L)$. En particulier, le radical de Baer est nilpotent.
\end{corollaire}

 Pour démontrer le théorème de Engel en toute généralité, nous devons d'abord préciser la notion d'automorphisme intérieur dans notre contexte. En effet, les dérivations ne sont pas des automorphismes; néanmoins, étant donné une dérivation nilpotente, on peut lui associer un automorphisme via la formule définissant l'exponentielle.
 \begin{fait}\label{exponentielle}
     Soit $L$ une algèbre de Lie sur un corps $K$ de caractéristique nulle. Alors :
     \begin{itemize}
         \item \cite{St} Toute dérivation $d$ nilpotente d'indice de nilpotence $n$ induit un automorphisme $\exp(d)=\sum_{i=0}^{n-1} \frac{d^i}{i!}$, d'inverse $exp(-d)$. En particulier, tout élément ad-nilpotent $x$ induit un automorphisme $\exp(\ad_x)=\exp(x)$.
         \item \cite{Hart} Soit $d$ une dérivation $k$-nilpotente, alors pour tout $x\in L$, il existe $a_1,\dots,a_n\in K$, tels que $d(x)=\sum_{n=1}^{k}a_nexp(nd)(x)$.
     \end{itemize}
 \end{fait}
 On fait le lien entre nos deux notions de normalisation. 
 \begin{Lemme}\label{normalisation}
     Soit $L$ une algèbre de Lie de caractéristique nulle. Soit $B$ un sous-groupe qui ne contient que des éléments ad-nilpotents (au sens de $L$) et soit $A$ un sous-espace. Alors $B\subseteq N_L(A)$ ssi $\exp(b)(A)=A^{\exp(b)}=A$ pour tout $b\in B$.
 \end{Lemme}
 \begin{proof}
     Soit $b\in N_L(A)$, $A$ est stable sous l'action de $\ad_b$ et de $\ad_{-b}$. Par conséquent, $\exp(b)$ et $\exp(-b)$ définissent sur $A$ des applications linéaires injectives inverses l'une de l'autre; ainsi, $A^{\exp(b)}=A$. Réciproquement, d'après le Fait \ref{exponentielle}, $\ad_b$ stablise $A$.
 \end{proof}
 \begin{definition}
     Soit $L$ une algèbre de Lie de caractéristique nulle. On appelle groupe d'automorphismes intérieurs, noté $\Int(L)$, le sous-groupe de $Aut(L)$ engendré par les éléments de la forme $\exp(x)$, où $x$ est un élément ad-nilpotent de $L$. Pour marquer l'analogie avec les groupes, étant donné un sous-espace $A$, on notera $A^f$ pour $f(A)$, où $f\in \Int(L)$.
 \end{definition}
 Désormais, nous sommes en position d'adapter de manière assez directe les arguments que Wagner donne pour les groupes $\mathfrak{M}_c$ dans \cite{Wag2}.
 \begin{Lemme}[{\cite[Lemme 1.8]{Wag2}}]\label{conjugaison nilpotence} 
    Soit $L$ une algèbre de Lie de caractérisque nulle et soit $S$ une sous-algèbre nilpotente de classe $c$. Soit $H$ une sous-algèbre qui ne contient que des éléments ad-nilpotents (au sens de $L$) et soit la sous-algèbre $N= \engendre{S^h : h\in \Int(H)}$. Si $H\leq \bigcap_{i\leq c}N_L(C_L^i(S))$, alors $N$ est nilpotente.
 \end{Lemme}
 \begin{proof}
     Soit $I=(x_0,x_1,\dots)$ une suite d'éléments dans $\bigcup_{\Int(H)} S^h$. Posons $y_0=x_0$ et $y_{i+1}=[y_i,x_{i+1}]$ pour $i\geq 0$. Remarquons que $C_L^i(S^h)=C_L^i(S)^h=C_L^i(S)$ (on obtient la dernière égalité en utilisant le Lemme \ref{normalisation}).

     Puisque $S$ est c-nilpotente, on a $S\leq C_L^c(S)$, et donc $y_0\in C_L^c(S)$ . De plus, si $y_i\in C_L^{c-i}(S)$ et $x_{i+1}\in S^h$ pour un certain $h\in \Int(H)$, alors $y_i\in C_L^{c-i}(S^h)$, ainsi $y_{i+1}\in C_L^{c-i-1}(S^h)=C_L^{c-i-1}(S)$ pour tout $i\leq c-1$. Par conséquent, $y_c=0$. Or, par anti-commutativité du crochet, tout commutateur de longueur $c+1$ à éléments dans $\bigcup_{\Int(H)} S^h$ est nul, et donc $N$ est nilpotente (car tout commutateur de $N$ est la somme de commutateurs à éléments dans $\bigcup_H S^h$).
 \end{proof}
 \begin{theoreme} [{à comparer avec \cite[Theorem 2.1]{Wag2}}]
     Soit $L$ une algèbre de Lie $\mathfrak{M}_c$  de caractéristique nulle. Supposons que $L$ est ad-nilpotente. Alors $L$ est nilpotente.
 \end{theoreme}
 \begin{proof}
 Supposons que $L$ n'est pas nilpotente. Remarquons que d'après le Théorème \ref{localement nilpotent} et le Corollaire \ref{ad-nilpotence résoluble}, une sous-algèbre est localement nilpotente ssi elle est résoluble ssi elle est nilpotente. Soit $S$ une sous-algèbre localement nilpotente maximale (qui existe d'après le Lemme de Zorn); il s'agit également d'une sous-algèbre résoluble maximale. Alors $N_L(S)=S$ : en effet, si $x\in N_L(S)$, alors la sous-algèbre $\engendre{S,x}=S+\engendre{x}$ est résoluble, donc égale à $S$ par maximalité.

     Il existe au moins une sous-algèbre nilpotente maximale $T$ distincte de $S$, car autrement, $S$ serait invariante sous l'action de $\Int(L)$. Par conséquent, d'après le Lemme \ref{normalisation}, $L=N_L(S)=S$, contradiction. 

     On considère 
     $I=\{ S\cap T : T$ est une sous-algèbre (localement) nilpotente maximale distincte de $S\}$.
     
    Soit $(K_i)$ une chaine ordonnée par l'inclusion et on note $K=\bigcup_i K_i$; c'est une sous-algèbre (localement) nilpotente de $S$, de classe $c$ disons. Soit $i_0$ tel que $C_L(K_{i_0}^j)$ soit minimal pour $j=1,\dots,c$ (et donc $C_L(K_i^j)=C_L(K_{i_0}^j)$ pour $i\geq i_0$). Soit $T$ une sous-algèbre nilpotente distincte de $S$ telle que $K_{i_0}=T\cap S$. D'après la condition du normalisateur pour les algèbres nilpotentes, il existe $y\in N_T(K_{i_0})-K_{i_0}$. On considère la sous-algèbre $F=\engendre{K_{i_0},y}=K_{i_0}+\engendre{y}$ qui est résoluble et normalise $K_{i_0}$ et $C_L^j(K_{i_0})$ pour $j\leq c$. Mais pour $j=1,\dots,c$, \[C_L(K^j)=C_L((\bigcup_i K_i)^j)=C_L\left(\bigcup_i K_i^j\right)=\bigcap_i C_L(K_i^j)=C_L(K_{i_0}^j).\] 

    Par conséquent, $C_L^j(K)=C_L^j(K_{i_0})$ pour tout $j\leq c$, d'après le Fait \ref{cent it 2}; ainsi, $F\leq \bigcap_{j\leq c} N_L(C_L^j(K))$.
    
    D'après le Lemme \ref{conjugaison nilpotence}, la sous-algèbre  $N=\engendre{K^f : f\in \Int(F)}$ est nilpotente (et contient $K$). Comme $F$ est résoluble et normalise $N$, on obtient que $\engendre{N,F}=N+F$ est une sous-algèbre résoluble contenant $K$ et $y$. Donc $\engendre{K,y}$ est résoluble (et localement nilpotente); elle est contenue dans une sous-algèbre localement nilpotente maximale $T_1\neq S$. Il suffit désormais de prendre $T_1\cap S$ pour majorant de la chaîne.

     Soit $I$ une intersection maximale de la forme $S\cap T$ (qui existe d'après le Lemme de Zorn); c'est une sous-algèbre nilpotente de $S$ et $T$. En vertu de la condition du normalisateur, on trouve $x\in N_S(I)-I$ et $y\in N_T(I)-I$. Il existe un élément $u$ de la forme $u=\ad_x^n(y) \in T$ qui n'est pas dans $T\cap S$ de longueur $n$ maximale (car $x$ est ad-nilpotent). Par conséquent, $\ad_x^{n+1}(u)=\ad_x(u)=-\ad_u(x)\in I=T\cap S$, et donc $\ad^k_u(x)\in I$ pour tout entier $k$ (puisque $u$ normalise $I$); en particulier, $x^{\exp(u)}\in S^{\exp(u)}\cap S$. Ainsi, on déduit du Lemme \ref{normalisation} que $S\cap S^{\exp(u)}>I$ et par maximalité $S=S^{\exp(u)}$. De même, on vérifie que pour tout $n\in \mathbb{Z}$, $S=S^{\exp(nu)}$. D'après le Lemme \ref{normalisation}, on obtient que $u\in N_L(S)=S$, contradiction.
 \end{proof}
 \begin{Remarque}
     Rosengarten démontre dans \cite{Ros} un analogue du théorème de Engel pour les anneaux de Lie de rang de Morlet fini (sans hypothèse sur la caractéristique). La démonstration repose sur l'existence d'une bonne notion de dimension à valeur entière pour les sous-ensembles définissables, et elle ne semble pas pouvoir se généraliser au cas $\omega$-stable et \textit{a fortiori} au cas $\mathfrak{M}_c$. Cependant, il resterait à élucider sa généralisation potentielle dans les anneaux de Lie fini-dimensionnels au sens de \cite{Wag3}.
 \end{Remarque}
 \section{Définissabilité dans le cas stable}
  On aborde maintenant la question de la définissabilité de l'idéal de Fitting dans le cas où l'algèbre de Lie est stable. Tout d'abord, nous adaptons à notre contexte et nous généralisons très légèrement \cite[Theorem 1.1.10]{Wag1}.
 \begin{proposition}\label{enveloppe}
     Soit $L$ une algèbre de Lie/un anneau de Lie stable. Alors tout idéal $N$ nilpotent (respectivement, résoluble) est contenu dans un sous-idéal définissable nilpotent (respectivement, résoluble).
 \end{proposition}
 \begin{proof}
    1) Par récurrence sur la classe de nilpotence, si $N$ est abélien, on considère $Z(C_L(N))$, qui est un idéal définissable abélien.
    \\
    On considère $Z(N)$, qui est contenu dans l'idéal définissable $A=Z(C_L(Z(N))$ (on remarque que $C_L(Z(N))$ contient $N$). On a que $N\leq C_L(A)$ et $Z(N)=A\cap N$; par conséquent, $N/Z(N)=N/A\cap N\simeq N+A/A$ est une sous-algèbre d'indice de nilpotence plus petit contenu dans $C_L(A)/A$. Par récurrence, il existe un idéal définissable $C/A$ de $C_L(A)/A$ contenant $N+A/A$, son image réciproque par le morphisme de projection est un idéal de $C_L(A)$ nilpotent définissable contenant $N$. Comme $C_L(A)$ est un idéal, il s'agit d'un sous-idéal.
    \\
    \\
    2) Si $N$ est résoluble, on prend le dernier terme de la série dérivée $A$ (qui est un idéal dans $L$) et on considère $B=Z(C_L(A))$ qui est un idéal définissable abélien contenant $A$. Or, $N\leq N_L(B)$ et donc $N+B/B\leq N_L(B)/B$ et on conclut par récurrence ($N_L(B)$ est un idéal de $L$). 
 \end{proof}
 \begin{corollaire}
     Soit $L$ une algèbre de Lie stable de caractéristique nulle; alors l'idéal de Fitting est définissable et nilpotent. 
 \end{corollaire}
 \begin{proof}
      On sait que l'idéal de Fitting est nilpotent, et il est donc contenu dans un sous-idéal définissable nilpotent. Mais d'après le Lemme \ref{sous-idéal nilpotent}, un tel sous-idéal est contenu dans un idéal nilpotent, et lui est donc égal par maximalité. 
 \end{proof}
 \begin{Remarque}
    Dans le cas d'un anneau de Lie $L$ de rang de Morley fini, le radical nilpotent et le radical résoluble existent et sont définissables. En effet, afin d'adapter la démonstration que donne Nesin pour les groupes dans \cite{Nes}, il suffit de vérifier que si $I$ est un idéal définissable de $L$, alors $I^{\circ}$ (la composante connexe au sens du groupe additif) est un idéal de $L$, et que, de plus, si $I$ est un idéal résoluble (respectivement, nilpotent), alors $I$ est contenu dans un idéal définissable résoluble (respectivement, nilpotent). Or, ces deux conditions sont bien satisfaites (l'argument suivant est dû à Adrien Deloro) : tout d'abord, étant donné un idéal définissable $I$ de $L$, on a, pour tout $x\in L$, que $\ad_x(I^{\circ})$ est un sous-groupe définissable connexe de $I$ et donc $\ad_x(I^{\circ})\leq_{gr} I^{\circ}$; de plus, si $I$ est un idéal de $L$, alors $d(I)$ (le plus petit sous-groupe définissable contenant $I$) est un idéal, car pour tout $x\in L$, $I\leq_{gr} \ad_x^{-1}(d(I))$ et donc $d(I)\leq_{gr}\ad_{x}^{-1}(d(I))$.\end{Remarque} 
 \printbibliography

 \end{document}